\documentclass[12pt]{article}
\usepackage{amsfonts,amssymb,amsthm}
\begin{document}
\newtheorem{thm}{Theorem}
\newtheorem{prop}{Proposition}
\newtheorem{plain}{Definition}
\newtheorem{rem}{Remark}
\newtheorem{cor}{Corollary}
\newtheorem{lemma}{Lemma}
\newtheorem{conj}{Conjecture}
\title{The $4^{th}$ structure.}
\date{}
\author{Yermolova -Magnusson J.
\footnote{e-mail: julia@oso.chalmers.se}\\
\small{\sl Onsala Space Observatory}\\
\small{\sl Chalmers University of Technology}\\
\small{\sl G\"oteborg, Sweden}\\
Stolin A. \footnote{e-mail: astolin@math.chalmers.se}\\
\small{\sl Department of Mathematics}\\
\small{\sl University of G\"oteborg}\\
\small{\sl G\"oteborg, Sweden}} \maketitle \abstract{In this paper we describe all Lie bialgebra structures on
the polynomial Lie algebra $\mathbf{g}[u]$, where $\mathbf{g}$ is a simple, finite dimensional, complex Lie
algebra. The results are based on an unpublished paper Montaner and Zelmanov  \cite{MZ}. Further, we introduce
quasi-rational solutions of the CYBE and describe all quasi-rational $r$-matrices for $\mathbf{sl}(2)$.}

\section{Introduction}

The aim of the present paper is to describe all Lie bialgebra structures on the polynomial Lie algebra
$P=\mathbf{g}[u]$, where $\mathbf{g}$ is a simple, finite dimensional Lie algebra over $\mathbb{C}$. The
starting point of the study became the following simple observation: up to classical twisting (defined in
\cite{KS}) there are only two Lie bialgebra structures on $\mathbf{g}$.  The corresponding Lie co-algebra
structures on these structures are given by $\delta_{1}=0$ and $\delta_{2}=[\Gamma_{DJ},a\otimes 1 +1\otimes
a]$, where $\Gamma_{DJ}$ is the classical Drinfeld-Jimbo modified $r$-matrix.

What is more difficult to prove is that of the classical twist does not change the corresponding classical
double, and for $\mathbf{g}$ we have, $\mathbf{D}_{1}(\mathbf{g})=\mathbf{g}\otimes \mathbb{C}[\varepsilon]$,
where $\varepsilon^{2}=0$, and $\mathbf{D}_{2}(\mathbf{g})=\mathbf{g}\otimes \mathbb{C}[\omega]$, where
$\omega^{2}=1$.(see \cite{S})

Similarly, according to an unpublished result of Montaner and Zelmanov, there are four Lie bialgebra structures
on $P=\mathbf{g}[u]$, again up to twisting. Therefore, we have to consider all the corresponding classical
doubles and then to derive Lie bialgebra structures related to one and the same classical double. So, it is
convenient to start with a description of the four possible cases for the double $\mathbf{D}(P)$ and the
corresponding Lie co-bracket on $P$.

\textbf{Case 1.}

In this case $\delta_{1}=0$ and $\mathbf{D}_{1}(P)=P+\varepsilon P^{*}$, where $\varepsilon^{2}=0$, and thus,
$\varepsilon P^{*}$ is a commutative ideal. The canonical symmetric nondegenerate invariant form $Q$ on
$\mathbf{D}_{1}(P)$ is given by the pairing between $P $ and $P^{*}$ :

$Q(p,q)=Q(p^{*},q^{*})=0$, $Q(p,\varepsilon q^{*})=q^{*}(p)$.

It is not difficult to show that there is a \textit{one-to-one} correspondence between Lie bialgebra structures
of the first type and finite-dimensional quasi-Frobenius Lie subalgebras of $P$.

\textbf{Case 2.}

Here $\delta_{2}(p(u))=[\Gamma_{2}(u,v),p(u)\otimes 1+1 \otimes p(v)]$ for $\Gamma_{2}(u,v)=\frac{\Omega}{u-v}$,
where $\Omega\in \mathbf{g}\otimes \mathbf{g}$ is the derived Casimir element and
$\mathbf{D}_{2}(P)=\mathbf{g}((u^{-1}))$. The canonical form $Q$ is given by the formula

$Q(p(u),q(u))=Res_{u=0}K(p,q)$,

where $K(p,q)$ is the Killing form of the Lie algebra $\mathbf{g}((u^{-1}))$ considered over
$\mathbb{C}((u^{-1}))$. In this case $P^{*}$ can be identified with $u^{-1}\mathbf{g}[[u^{-1}]]$. The
corresponding Lie bialgebra structures are in a \textit{one-to-one} correspondence with \textit{rational
solutions} of the CYBE (see \cite{S1}).

\textbf{Case 3.}

In the third case the $1-cocycle$ $\delta_{3}$ is defined as $\delta_{3}(p(u))=[\Gamma_{3}(u,v),p(u)\otimes 1+1
\otimes p(v)]$ with the corresponding  $\Gamma_{3}(u,v)=\frac{v\Omega}{v-u}+r_{DJ}$, where $r_{DJ}$ is the
classical Drinfeld-Jimbo modified $r$-matrix. The classical double has the following form
$\mathbf{D}_{3}(P)=\mathbf{g}((u^{-1}))\oplus \mathbf{g}$. The expression for the invariant bilinear $Q$ is:

$Q((f(u),a),(g(u),b))=Res_{u=0}u^{-1}K(f,g)-K(a,b)$, where $a,b \in \mathbf{g}$

The Lie algebra $P^{\ast}$ is identified with:

$u^{-1}\mathbf{g}[[u^{-1}]]\oplus \{(l,k)\in \mathbf{b}_{+}\oplus \mathbf{b}_{-}:
l_{\textbf{h}}+k_{\textbf{h}}=0\} $,

where $\textbf{h}$ is a fixed Cartan subalgebra, $\mathbf{b}_{\pm}$ are the corresponding Borel subalgebras of
$\mathbf{g}$ and $l_{\textbf{h}}$, $k_{\textbf{h}}$ is a Cartan part of $l$ and $k$ respectively. There is a
natural \textit{one-to-one} correspondence between these Lie bialgebra structures with the so-called
\textit{quasi-trigonometric solutions}  of the CYBE. \cite{Y}

\textbf{Case 4.}

This is a new case. The Lie bialgebra structure on $P$ is given by:

$\delta_{4}(p(u))=[\Gamma_{4}(u,v),p(u)\otimes 1+1 \otimes p(v)]$, with
$\Gamma_{4}(u,v)=\frac{\Omega}{u^{-1}-v^{-1}}=\frac{uv\Omega}{v-u}$.

$\mathbf{D}_{4}(P)=\mathbf{g}((u^{-1}))\oplus \mathbf{D}_{1}(\mathbf{g})$ , where $\mathbf{D}_{1}(\mathbf{g})$
was defined above.

Let us describe the corresponding form $Q$:

$Q(f(u)+A_{0}+A_{1}\varepsilon,g(u)+B_{0}+B_{1}\varepsilon)$,

where $f(u)=\sum_{-\infty}^{N}a_{k}u^{k}$, $g(u)=\sum^{M}_{-\infty}b_{k}u^{k}$ and $a_{i},b_{i},A_{i},B_{i}\in
\mathbf{g}$.
\begin{equation}
Q(f(u)+A_{0}+A_{1}\varepsilon,g(u)+B_{0}+B_{1}\varepsilon)=\mathbf{Res}_{u=0}u^{-2}K(f,g)-K(A_{0},B_{1})-K(A_{1},B_{0}).
\end{equation}
$P$ is embedded into $\mathbf{D}_{4}$ as follows:
$i(\sum_{k=0}^{N}c_{k}u^{k})=\sum_{k=0}^{N}c_{k}u^{k}+c_{0}+c_{1}\varepsilon$ (It is easy to check that $i(P)$
is a Lagrangian subalgebra of $\mathbf{D}_{4}(P)$ )

$P^{*}$ can be identified with $\mathbf{g}[[u^{-1}]]\oplus \mathbf{g}\varepsilon.$ In order to verify  that
$\mathbf{D}_{4}(P)$ corresponds to $\Gamma_{4}(u,v)=\frac{uv\Omega}{v-u}$, we have to find the dual bases in $P$
and $P^{*}$ with respect to $Q$ and compute the corresponding dual base element.

Let us do this.

Let $\{I_{i}\}$ be an orthonormal basis of $\mathbf{g}$ with respect to $K$. Then $i(\mathbf{g}[u]) \subset
\mathbf{D}_{4}$ poses the following basis: $\{I_{i}u^{k} (k\geq 2),(I_{i}u, I_{i}\varepsilon), (I_{i},I_{i})\}$.
It is not difficult to find the dual basis in $P^{\ast}$: $\{I_{i}u^{-k+1}(k\geq 2), (I_{i},0),(0,-\varepsilon
I_{i})\}$. The corresponding dual base element is:

\begin{equation}
\sum_{i} \sum_{k=2}^{\infty}I_{i}u^{k}\otimes
I_{i}v^{-k+1}+\sum_{i}(I_{i}u,I_{i}\varepsilon)\otimes(I_{i},0)+\sum_{i}(I_{i},I_{i})\otimes(0,-\varepsilon
I_{i}).
\end{equation}

If we project this element onto the first component in the decomposition

$\mathbf{D}_{4}(P)=\mathbf{g}((u^{-1}))\oplus \mathbf{D}_{1}(\mathbf{g})$, we get:

\begin{equation}
\sum_{i} \sum_{k=1}^{\infty}I_{i}u^{k}\otimes I_{i}v^{-k+1}=u\sum_{i}\sum_{m=0}^{\infty}I_{i}u^{m}\otimes
I_{i}v^{-m}=\frac{uv\Omega}{v-u}.
\end{equation}

Therefore, we have proved the following result:
\begin{thm}
The double $\mathbf{D}_{4}(P)$ corresponds to the co-bracket on $P$ given by the formula
\begin{equation}
\delta_{4}(p(u))=[\Gamma_{4}(u,v),p(u)\otimes 1+1 \otimes p(v)],
\end{equation}
with $\Gamma_{4}(u,v)=\frac{uv\Omega}{v-u}$
\end{thm}

\begin{plain}
Solutions of the classical Yang-Baxter equation of the form $q(u,v)=\frac{uv\Omega}{v-u}+p(u,v)$, where $p(u,v)$
is a skew-symmetric polynomial, are called \textbf{quasi-rational} $r$- matrices.
\end{plain}

\section{Quasi-rational solutions of CYBE.}

Following methods developed in \cite{KPST},\cite{KS},\cite{S}, we can prove that:

\begin{thm}

(i). All quasi-rational $r$-matrices provide one and the same double $\mathbf{D}_{4}(P).$

(ii). Quasi-rational solutions of the CYBE are in a one-to-one correspondence with Lagrangian subalgebras
$\mathbf{W}\subset \mathbf{D}_{4}(P)$, satisfying the following conditions:

1. $\mathbf{W}\cap P= \{0\}$

2. $\mathbf{W}\oplus P=\mathbf{D}_{4}(P)$

3. $\mathbf{W}\supset u^{-N}\mathbf{g}[[u^{-1}]]$, for some $N>0$.
\end{thm}
Since $P\subset \mathbf{D}_{4}(P)$, we can treat the group $\mathbf{Ad}(P)$, as a subgroup of
$\mathbf{Ad}(\mathbf{D}_{4}(P))$. If $q(u,v)$ is a quasi-rational $r$-matrix, then $\mathbf{Ad}(p(u)\otimes
p(v)q(u,v))$ is also a quasi-rational solution of the CYBE.

We say that two quasi-rational $r$-matrices $q_{1}(u,v)$ and $q_{2}(u,v)$ are gauge equivalent if there exists
$p(u)\in \mathbf{Ad}\mathbf{g}[u]$ such that:

$q_{2}(u,v)=\mathbf{Ad}(p(u)\otimes p(v))q_{1}(u,v)$.

It is not difficult to prove that:

\begin{prop}
Quasi-rational $r$-matrices $q_{1}(u,v)$ and $q_{2}(u,v)$ are gauge equivalent iff
 $\mathbf{W}_{2}=\mathbf{Ad}(p(u))\mathbf{W}_{1}$, where $\mathbf{W}_{i}$ corresponds to $q_{i}(u,v)$,
$(i=1,2)$.
\end{prop}
\begin{thm}
Let $\mathbf{g}=\mathbf{sl}(n)$. For any quasi-rational $r$-matrix $q_{1}(u,v)$, there exists another
quasi-rational $r$-matrix $q_{2}(u,v)$, which is gauge equivalent to $q_{1}(u,v)$, and such that the
corresponding $\mathbf{W}_{2}\subset \mathbf{D}_{4}(P)$ is contained in $d_{k}^{-1}\mathbf{sl}(n,
\mathbb{C}[[u^{-1}]])d_{k}\oplus \mathbf{sl}(n,\mathbb{C}[\varepsilon])$.

Here $d_{k}=diag (\underbrace{1,...,1}_{k},\underbrace{u,...,u}_{n-k})$.
\end{thm}

Now we are ready to give a more explicit description of the quasi-rational $r$ -matrices.
\begin{prop}
If $\mathbf{W}$ corresponds to a quasi-rational $r$-matrix $q(u,v)$ and $\textbf{W}\subset
\mathbf{sl}(n,\mathbb{C}[[u^{-1}]])\oplus \mathbf{sl}(n, \mathbb{C}[\varepsilon])$, then $q(u,v)=\frac{uv\Omega
}{v-u}+q$, where $q$ is a skew-symmetric solution of the CYBE. Conversely, if $q(u,v)=\frac{uv\Omega }{v-u}+q$
is a quasi-rational $r$-matrix, then the corresponding $\mathbf{W}$ is contained in
$\mathbf{sl}(n,\mathbb{C}[[u^{-1}]])\oplus \mathbf{sl}(n, \mathbb{C}[\varepsilon])$.
\end{prop}

\begin{rem}
Both statements above are valid for any $\mathbf{g}$.
\end{rem}
\begin{rem}
Since skew-symmetric $r$-matrices are in a one-to-one correspondence with quasi-Frobenius subalgebras of
$\mathbf{g}$, we have proved that there is a one-to-one correspondence between quasi-rational $r$-matrices of
the form $\frac{uv\Omega }{v-u}+q$ and quasi-Frobenius subalgebras of $\mathbf{g}$.
\end{rem}

We continue with the following:

\begin{lemma}
Let $\mathbf{W}_{k}=d_{k}^{-1}\mathbf{sl}(n, \mathbb{C}[[u^{-1}]])d_{k}\oplus
\mathbf{sl}(n,\mathbb{C}[\varepsilon])$. Then its orthogonal complement $\mathbf{W}_{k}^{\perp}$ with respect to
the form $Q$ on $\mathbf{D}_{4}(\mathbf{sl}(n))$,  is isomorphic to

$d_{k}^{-1}\mathbf{sl}(n, \mathbb{C}[[u^{-1}]])d_{k}$ and $\frac{\mathbf{W}_{k}}{\mathbf{W}_{k}^{\perp}}\cong
\mathbf{sl}(n,\mathbb{C}[\varepsilon])$.
\end{lemma}

Now we can to describe all the quasi-rational $r$-matrices related to $\mathbf{W}_{k}$.

Let $\mathbf{W}$ correspond to $q(u,v)$ and let $\mathbf{W}\subset\mathbf{W}_{k}$. Since $\mathbf{W}$ is
Lagrangian, we get $\mathbf{W}_{k}^{\perp}\subset \mathbf{W}\subset \mathbf{W}_{k}$. Therefore, the image of
$\mathbf{W}$ in $\frac{\mathbf{W}_{k}}{\mathbf{W}_{k}^{\perp}}$ is a Lagrangian subalgebra , which is
transversal to $P\cap \mathbf{W}_{k}$ (it is easy to see that the canonical map $P\cap
\mathbf{W}_{k}\longrightarrow \frac{\mathbf{W}_{k}}{\mathbf{W}_{k}^{\perp}} $ is an embedding, since $P\cap
\mathbf{W}_{k}=\{0\}$ and $\mathbf{W}\supset\mathbf{W}_{k}^{\perp}$). The image of $P\cap \mathbf{W}_{k}$ in
$\frac{\mathbf{W}_{k}}{\mathbf{W}_{k}^{\perp}}\cong \mathbf{sl}(n,\mathbb{C}[\varepsilon])$ can be computed.

Let $P_{k}$ be the maximal parabolic subalgebra of $\mathbf{sl}(n)$, which contains $B^{+}$ and corresponds to
the maximal root $\alpha_{k}$. Then the image of $P\cap \mathbf{W}_{k}$ is isomorphic to $P_{k}+\varepsilon
P_{k}^{\perp}$, where $P_{k}^{\perp}$ is the orthogonal complement to $P_{k}$ with respect to the Killing form
on $\mathbf{sl}(n)$. Thus, we proved:

\begin{thm}
There is a one-to-one correspondence between quasi-rational $r$-matrices related to $\mathbf{W}_{k}$ and
Lagrangian subalgebras of $\mathbf{sl}(n, \mathbb{C}[\varepsilon])$ transversal to $P_{k}+\varepsilon
P_{k}^{\perp}$.
\end{thm}

At this point we note that \textit{Theorem 4} together with results of \cite{S}, provide a \textit{one-to-one}
correspondence between rational and quasi-rational $r$-matrices in $\mathbf{sl}(n)$.

Therefore, the following result is proved.

\begin{thm}
Quasi-rational $r$-matrices related to $\mathbf{W}_{k}$ are in a one-to-one correspondence with the pairs
$(L,B)$, where $L\subset\mathbf{sl}(n)$ is a subalgebra such that $L+P_{k}=\mathbf{sl}(n)$, and $B$ is a
2-cocycle on $L$, such that $B$ is non-degenerate on $L\cap P_{k}$. (see \cite{S} for details).
\end{thm}

Finally, we would like to compute all quasi-rational $r$-matrices for $\mathbf{sl}(2)$. In this case $k=0,1$.

If $k=0$, then the corresponding quasi-rational $r$-matrices are in a \textit{one-to-one} correspondence with
the quasi-Frobenius subalgebras of $\mathbf{sl}(2)$. We have two cases:

1. $L=0\Longrightarrow q_{0}(u,v)=\frac{uv\Omega}{v-u}$

2. $L=B^{+}\Longrightarrow q_{1}(u,v)=\frac{uv\Omega}{v-u}+e\otimes h-h\otimes e$

If $k=1$, then we have only one pair $(L,B)$ satisfying conditions $L+P_{1}=\mathbf{sl}(2)$, $B$ is
non-degenerate on $L\cap P_{1}$: $L=\mathbf{sl}(2)$, $B(x,y)=K(f,[x,y])$. Here, $\{e,f,h\}$ is the standard
basis of $\mathbf{sl}(2)$ and $K$ is the Killing form on $\mathbf{sl}(2)$ .

The rational $r$-matrix corresponding to this data is:
\begin{equation}
\frac{\Omega}{u-v}+eu\otimes h-h\otimes ev.
\end{equation}
However, the quasi-rational $r$-matrix, which corresponds to the same data is rather different:

\begin{equation}
q_{2}(u,v)=\frac{uv\Omega}{v-u}+\frac{1}{2}h\otimes e-\frac{1}{2}e\otimes h-eu\otimes f+f\otimes ev.
\end{equation}

One can prove that up to gauge equivalence $q_{0}(u,v)$, $q_{1}(u,v)$ and $q_{2}(u,v)$ exhaust all
quasi-rational $r$-matrices for $\mathbf{sl}(2)$.

\bigskip
{\small \textbf{Acknowledgements}: Part of this paper was written during the visit of the first author to
Department of Mathematics and Statistics, University of Sao Paulo. The author is very grateful for hospitality
and useful conversations. }

\end{document}